\documentclass[conference]{IEEEtran}
\IEEEoverridecommandlockouts
\usepackage{cite}
\usepackage{amsmath,amssymb,amsfonts}
\usepackage{algorithmic}
\usepackage{graphicx}
\usepackage{textcomp}
\usepackage{xcolor}
\usepackage{mathtools}
\usepackage{float}
\newcommand{\mr}{\mathrm}                   

\newcommand{\veg}[1]{\bm{#1}}               
\newcommand{\mat}[1]{\mathsfbfit{#1}}       
\newcommand{\uv}[1]{\hat{\veg{#1}}}         
\newcommand{\op}[1]{\mathcal{#1}}           
\newcommand{\matel}[1]{\begin{bmatrix} #1 \end{bmatrix}}    

\newcommand{\matI}{\mathbfsf{I}}
\newcommand{\matO}{\mathbfsf{0}}

\newcommand{\dd}{\mathrm{d}}  

\newcommand{\jm}{\mathrm{j}}  

\newcommand{\e}{\mathrm{e}}

\newcommand{\T}{\mr{T}}

\usepackage[utf8]{inputenc}

\usepackage[T1]{fontenc}
\usepackage{subcaption}
\usepackage{float}
\usepackage{tikz}
\usepackage{pgfplots}
\usepackage{pgfplotstable}
\usepackage{siunitx}
\usepackage{cleveref}
\pgfplotsset{compat=newest}

\usepackage{amsthm}
\DeclareMathAlphabet{\mathbfsf}{\encodingdefault}{\sfdefault}{bx}{n}
\usepackage{bm}
\usepackage[OMLmathsfit]{isomath}

\newcommand{\vt}[1]{\boldsymbol{#1}}

\def\BibTeX{{\rm B\kern-.05em{\sc i\kern-.025em b}\kern-.08em
    T\kern-.1667em\lower.7ex\hbox{E}\kern-.125emX}}
\begin{document}

\title{A New Preconditioner for the EFIE Based on Primal and Dual Graph Laplacian Spectral Filters
%
%
}
%
%
\author{\IEEEauthorblockN{ Lyes Rahmouni* and Francesco P. Andriulli}
\IEEEauthorblockA{Department of Electronics and Telecommunications, Politecnico di Torino, Turin, Italy \\
lyes.rahmouni@polito.it, francesco.andriulli@polito.it}
}

\maketitle

\begin{abstract}
The Electric Field Integral Equation (EFIE) is notorious for its ill-conditioning both in frequency and h-refinement. Several techniques exist for fixing the equation conditioning problems based on hierarchical strategies, Calderon techniques, and related technologies. This work leverages on a new approach, based on the construction of tailored spectral filters for the EFIE components which allow the block renormalization of the EFIE spectrum resulting in a provably constant condition number for the equation. This is achieved without the need for a barycentric refinement and with low computational overhead compared with other schemes. In particular, only sparse matrices are required in addition to the EFIE original matrix. Numerical results will show the robustness of our scheme and its application to the solution of realistic problems.
\end{abstract}
%
%
\begin{IEEEkeywords}
EFIE, preconditioning, projectors, loop-star decomposition
\end{IEEEkeywords}

%
%

\section{Introduction}

Numerous problems in electromagnetics can be formulated as boundary integral equations. Surface integral equations are especially efficient because they require the discretization of scatterers boundaries only but, however, they give rise to dense linear systems that are generally solved with iterative methods. The quality of the solution as well as the number of iterations required to achieve convergence intrinsically depends on the spectral properties of the discretized operators. Unfortunately, a large subset of integral formulations produce matrices whose condition number grows unbounded with respect to the inverse of the mesh parameter $h$ (the average edge length) \cite{andriulli2008multiplicative}. This ill-conditioning is mainly due to the fact that some of the underlying continuous operators are Fredholm integral operators of the first-kind, which are known to have eigenvalues accumulating at zero and/or at infinity \cite{andriulli2010solving,andriulli2009hybrid}. This is a critical limitation as it becomes prohibitively expensive to solve problems of practical interest. At the same time, electromagnetic integral equations are often unstable when the frequency decreases. This is mainly due to an unfavorable scaling of the operator components and to a numerical cancellation phenomenon of the electric current solutions. 

Several efforts have been profused in developing strategies for addressing the above problems. On the one hand algebraic preconditioners have been proposed (see for example \cite{lee2003incomplete} and \cite{benzi1996sparse}) which, although improving the conditioning properties, they are still showing growing condition numbers for  finer meshes. A second class of preconditioners are tuned to the spectral properties of the involved operators and encompass hierarchical, Calderon strategies, and related methods (\cite{andriulli2008multiplicative, adrian2019refinement} and references therein). These schemes, however, may often results in computational overheads either because a barycentric refinement is required or because extra dense operators must be computed.

A different approach is adopted here: we leverage on the design suitably conceived spectral filters for the EFIE components which allow for the block renormalization of the EFIE spectrum. This results in a constant condition number without the need for a barycentric refinement and with low computational overhead compared with other schemes. In particular, only sparse matrices are required in addition to the EFIE original matrix. Numerical results will show the robustness of our scheme and its application to the solution of realistic problems.

%
%
\section{Background and Notation}

 Let $\Gamma$ be a Lipschitz boundary representing the surface of a Perfect Electrically Conducting (PEC) object and $\uv{n}$ its outward pointing unit normal. A time harmonic incident wave $\vt{E}^i$ induces a surface electric current density $\vt{J}$, which in turn generates a scattered field $\vt{E}^s$. The latter can be computed by solving the EFIE which reads
 \begin{equation}
  \label{EFIE}
  	-\uv{n}(\vt{r}) \times \vt{E^i}=\op{T} \vt{J}=\op{T}_A \vt{J} +  \op{T}_\phi \vt{J}
  \end{equation}
  where the $\op{T}_A \vt{J}$ is the vector potential
\begin{equation}
\op{T}_A \vt{J} = \uv{n}(\vt{r}) \times  \jm k  \int_\Gamma \frac{\e^{\jm k\| \vt{r} - \vt{r'} \|}}{4\pi \|\vt{r} - \vt{r'} \|}   \vt{J}(\vt{r'})  \dd S(\vt{r'}) 
\end{equation}
and $\op{T}_\phi \vt{J}$ is the scalar potential
\begin{equation}
		 \op{T}_\phi \vt{J} = -\uv{n}(\vt{r}) \times \frac{1}{\jm k} \nabla_{\vt{r}} \int_\Gamma \frac{\e^{\jm k\| \vt{r} - \vt{r'} \|}}{4\pi\|\vt{r} - \vt{r'} \|}  \nabla_{\vt{r'}} \cdot \vt{J}(\vt{r'})  \dd S(\vt{r'}) 
\end{equation}
in which $k=\omega \sqrt{\epsilon\mu}$ is the wavenumber. Following standard discretization strategy, we approximate $\Gamma$ with a triangular mesh of average edge length $h$. On this mesh, the current density $ \vt{J}$ is expanded  as $ \vt{J} = \sum\nolimits_{n = 1}^N {I_n } \vt{f}_n (\vt{r}) $, where $\vt{f}_n (\vt{r})$ are $N$ Rao-Wilton-Glisson (RWG) basis functions \cite{rao1982electromagnetic}. Equation (\ref{EFIE}) is then tested with $\uv{n}(\vt{r}) \times \vt{f}_n (\vt{r})$ to obtain a linear system $\mat{T} \mat{j} =\mat{e}$ where $\mat{T} = \jm k \mat T_A + (\jm k)^{-1} \mat T_\phi$, with $\mat{T}_A = \left\langle {\uv{n}(\vt{r}) \times \vt{f}_n (\vt{r}), \op{T}_A(\vt{f}_n (\vt{r}))} \right\rangle $, $\mat{T}_\phi = \left\langle {\uv{n}(\vt{r}) \times \vt{f}_n (\vt{r}), \op{T}_\phi(\vt{f}_n (\vt{r}))} \right\rangle $ and $\mat{e}=\left\langle{ \vt{f}_n (\vt{r}), -\vt{E}^i} \right\rangle$. Unfortunately, this system is ill-conditioned both in low frequencies and with dense discretization, as $\text{cond} (\mat{T}) \lesssim 1/(hk)^2$.

%
%

Analyzing equation \eqref{EFIE}, it is clear that the two components of the operator $\op{T}$ scale inversely in frequency, which is the source of the low frequency breakdown of the EFIE. Using a loop-star decomposition \cite{ andriulli2012loop}, it is possible to decouple the contributions of the solenoidal and irrotational  currents, which will allow for a diagonal preconditioner to effectively cure the low frequency breakdown. In particular, multiplying left and right with $\mat{\Upsilon}= [\mat{\Lambda}/\sqrt{k} \quad \mat{\Sigma}\sqrt{k}]$, where $\mat{\Lambda}$ and $\mat{\Sigma}$ are the RWG to loop and RWG to star transformation matrices, eliminates the unfavorable frequency ill-scaling. The resulting matrix $\mat{T}_{LS}=\mat{\Upsilon}^\text{T} \mat{T} \mat{\Upsilon}$ has a condition number stable with decreasing frequency. This, however, comes at the expense of a degraded conditioning in dense discretization. Indeed, the condition number of the new EFIE system grows cubically with the edge parameter, that is,  $\text{cond}(\mat{T}_{LS}) \lesssim 1/h^3$. In the following, we present a new strategy to cure the dense discretization breakdown. 
%
%

\section{ A New Preconditioner Based on Primal and Dual Graph Laplacian Spectral Filters} 
Following the above-mentioned considerations, it is necessary to further regularize the operator $\mat{T_{LS}}$  $h$-dependency. The strategy adopted in this work consists in partitioning the spectrum of the operator $\mat{T_{LS}}$ into $N$ sub-intervals by following the natural ordering of the primal and dual graph Laplacian. 
Subsequently, we normalize each sub-interval eigenvalues by the largest eigenvalue of the $i^{th}$ interval $\sigma_{\max}^i$. A naive eigen-decomposition of the different operators, however, is not practical given its $\mathcal{O}(N^3)$ and $\mathcal{O}(N^2)$ complexity for dense and sparse matrices, respectively. A computationally efficient approach, proposed in this work, leverages on the sparsity of the graph Laplacian to build spectral filters $\mat Q_i$. We first introduce a family of low-pass spectral filters of length $2^l$  
\begin{align}
    \mat P_i^\Sigma &= \frac{\matI}{{\matI + \mat \upDelta^\Sigma_i}}  \quad \text{for}\,\, i=1 \dots N_{\Sigma}  \,,\\
    \mat P_i^\Lambda &= \frac{\matI}{{\matI + \mat \upDelta^\Lambda_i}}  \quad \text{for}\,\, i=1 \dots N_{\Lambda}  \,,
\end{align}
where $N_{\Sigma} = \log_2(\sigma_{\max}^\Sigma)$,  $N_{\Lambda} = \log_2(\sigma_{\max}^\Lambda)$ and
\begin{equation}
\mat \upDelta^{\Sigma, \Lambda}_i = \left( {\frac{\mat \upDelta^{\Sigma, \Lambda}}{2^i}}\right)^n
\end{equation}
in which the parameter $n$ controls the sharpness of the filters. The primal and dual Laplacian used in the filters are given by 
 \begin{equation}
 \mat{\upDelta^\Sigma} = \mat{\Sigma}^\T\mat{\Sigma}
 \end{equation}

 \begin{equation}
  \mat{\upDelta^\Lambda} = \mat{\Lambda}^\T  \mat{\Lambda}
 \end{equation}

The filters $\mat Q_i$ are then defined as 
\begin{equation}
\mat Q_i = (\mat P_i - \mat P_{i-1}) \frac{1}{\sqrt{\sigma_{\max}^i}}.
\end{equation}
The well-conditioned EFIE we propose is finally given by 
\begin{equation}
\matel{ \mat Q_i^\Sigma & \matO\\
        \matO   & \mat Q_i^\Lambda}
\matel{
k\widetilde{ \mat{\Sigma}}^\T(\mat{T_A}+\mat{T_\phi})\widetilde{ \mat{\Sigma}} & \widetilde{ \mat{\Sigma}}^\T \mat{T_A} \mat{\Lambda} \\
\mat{\Lambda}^\T\mat{T_A}\widetilde{ \mat{\Sigma}} & \frac{1}{ k}\mat{\Lambda}^\T \mat{T_A} \mat{\Lambda}
}
\matel{ \mat Q_i^\Sigma & \matO\\
        \matO   & \mat Q_i^\Lambda}
\end{equation}
where
\begin{equation}
   \widetilde{ \mat{\Sigma}}= \mat{\Sigma} \left( {\mat{\Sigma}^\T \mat{\Sigma}} \right)^{+}\
\end{equation}
The introduction of further Gram matrices can be used to further reduce  the condition number especially in case of non uniform meshes,  but we omitted this analysis here  for the sake of brevity. 

\section{Numerical Results}

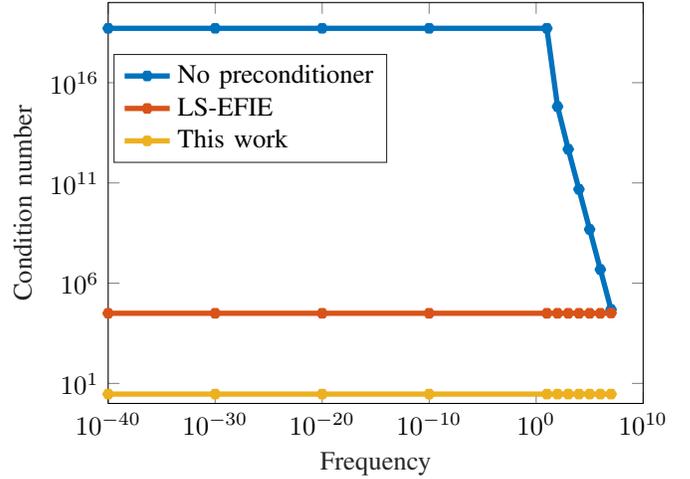
\begin{figure}
\centering
%
%
\definecolor{mycolor1}{rgb}{0.00000,0.44700,0.74100}%
\definecolor{mycolor2}{rgb}{0.85000,0.32500,0.09800}%
\definecolor{mycolor3}{rgb}{0.92900,0.69400,0.12500}%
\begin{tikzpicture}

\begin{axis}[%
width=2.8in,
height=2.1in,
at={(0.758in,0.481in)},
scale only axis,
xmode=log,
xmin=1e-40,
xmax=10000000000,
xminorticks=true,
xlabel style={font=\color{white!15!black}},
xlabel={Frequency},
ymode=log,
ymin=1,
ymax=1e+20,
yminorticks=true,
ylabel style={font=\color{white!15!black}},
ylabel={Condition number},
axis background/.style={fill=white},
xminorgrids,
yminorgrids,
legend style={at={(0.01,0.6)}, anchor=south west, legend cell align=left, align=left, draw=white!15!black}
]
\addplot [color=mycolor1, line width=2.0pt, mark=asterisk, mark options={solid, mycolor1}]
  table[row sep=crcr]{%
10000000	48169\\
1000000	4817000\\
100000	481700000\\
10000	48170000000\\
1000	4817400000000\\
100	    6.4632e+14\\
10	5.1705e+18\\
1e-10	5.1705e+18\\
1e-20	5.1705e+18\\
1e-30	5.1705e+18\\
1e-40	5.1705e+18\\
};
\addlegendentry{No preconditioner}

\addplot [color=mycolor2, line width=2.0pt, mark=asterisk, mark options={solid, mycolor2}]
  table[row sep=crcr]{%
10000000	32430\\
1000000	31485\\
100000	31475\\
10000	31475\\
1000	31475\\
100	31475\\
10	31475\\
1e-10	31475\\
1e-20	31475\\
1e-30	31475\\
1e-40	31475\\
};
\addlegendentry{LS-EFIE}

\addplot [color=mycolor3, line width=2.0pt, mark=asterisk, mark options={solid, mycolor3}]
  table[row sep=crcr]{%
10000000	2.9295\\
1000000	2.9294\\
100000	2.9294\\
10000	2.9294\\
1000	2.9294\\
100	2.9294\\
10	2.9294\\
1e-10	2.9294\\
1e-20	2.9294\\
1e-30	2.9294\\
1e-40	2.9294\\
};
\addlegendentry{This work}

\end{axis}
\end{tikzpicture}%
    \caption{Condition number as a function of the frequency}
    \label{Sph_Cond_VS_Freq}
\end{figure}

Our first numerical result is intended to demonstrate that the presented technique delivers an EFIE immune from the low frequency breakdown. To that end, we simulated a spherical PEC geometry of radius 1m, discretized with 3270 triangles. Plane waves of different frequencies were used as excitations. The condition numbers of $\mat{T}$, $\mat{T_{LS}}$ and our technique are reported in  \Cref{Sph_Cond_VS_Freq}. We can see that even though the Loop-Star decomposition solves the low frequency breakdown, our scheme further improves the condition number of the EFIE.

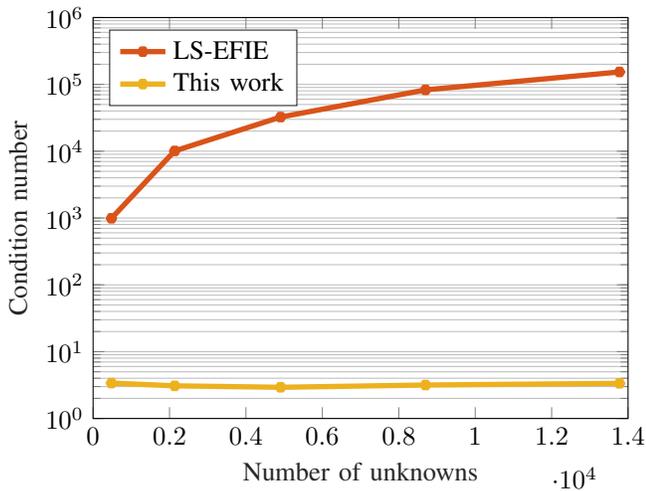
\begin{figure}
\centering
%
%
\definecolor{mycolor1}{rgb}{0.00000,0.44700,0.74100}%
\definecolor{mycolor2}{rgb}{0.85000,0.32500,0.09800}%
\definecolor{mycolor3}{rgb}{0.92900,0.69400,0.12500}%

\begin{tikzpicture}

\begin{axis}[%
width=2.8in,
height=2.1in,
at={(0.758in,0.481in)},
scale only axis,
xmin=0,
xmax=14000,
xlabel style={font=\color{white!15!black}},
xlabel={Number of unknowns},
ymode=log,
ymin=1,
ymax=1000000,
yminorticks=true,
ylabel style={font=\color{white!15!black}},
ylabel={Condition number},
axis background/.style={fill=white},
xminorgrids,
yminorgrids,
legend style={at={(0.03,0.97)}, anchor=north west, legend cell align=left, align=left, draw=white!15!black}
]
\addplot [color=mycolor2, line width=2.0pt, mark=asterisk, mark options={solid, mycolor2}]
  table[row sep=crcr]{%
480	990.562\\
2136	10089\\
4908	32430\\
8700	82694\\
13776	153400\\
};
\addlegendentry{LS-EFIE}

\addplot [color=mycolor3, line width=2.0pt, mark=asterisk, mark options={solid, mycolor3}]
  table[row sep=crcr]{%
480	3.3879\\
2136	3.0797\\
4908	2.9295\\
8700	3.1656\\
13776	3.353\\
};
\addlegendentry{This work}

\end{axis}
\end{tikzpicture}%
    \caption{Condition number as a function of number of elements}
    \label{Sph_Cond_VS_h}
\end{figure}

As second example, we considered a PEC sphere illuminated by a plane
wave oscillating at \SI{1}{\hertz}. \Cref{Sph_Cond_VS_h} shows the condition number of the EFIE matrix of the standard loop-tree preconditioner and our new scheme. We can see that the condition number of $\mat{T}_{LS}$ grows unbounded as the discretization increases. Our formulation offers a regularized and stable solver.

\begin{figure}
	\centering
	\includegraphics[trim = 0mm 0 0 0mm, clip, width=0.4\textwidth]{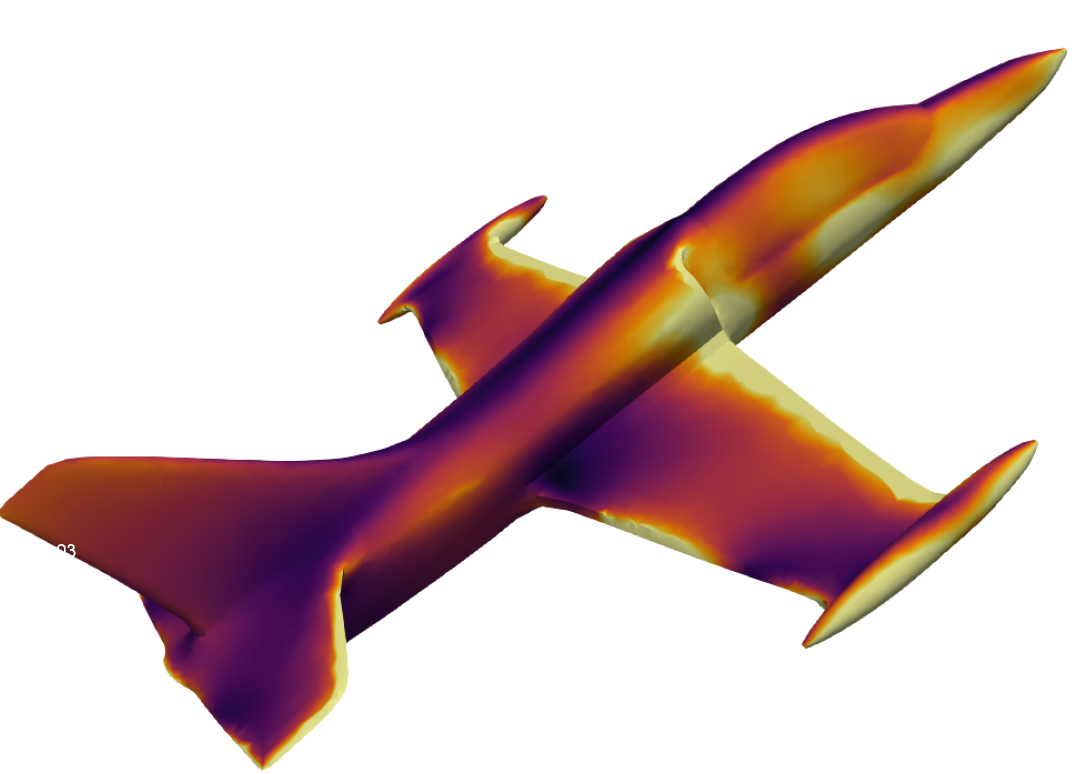}
	\caption{The induced surface current}
	\label{fig:Fighter}
\end{figure}

 The last example shows the study of a real case scenario. Accordingly, we simulated the aircraft shown in \cref{fig:Fighter}, discretized into $16400$ triangles. As an excitation, we considered a plane wave oscillating at \SI{e5}{Hz}. Using loop-star decomposition, the conjugate gradient converged in 2400 iterations, while with our techniques, it took 330 iterations.

\section*{Acknowledgment}

This work was supported by the European Research Council (ERC) under the European Union’s Horizon 2020 research and innovation program (grant agreement No 724846, project 321).

\bibliographystyle{IEEEtran}
\bibliography{references}

\begin{thebibliography}{1}
\providecommand{\url}[1]{#1}
\csname url@samestyle\endcsname
\providecommand{\newblock}{\relax}
\providecommand{\bibinfo}[2]{#2}
\providecommand{\BIBentrySTDinterwordspacing}{\spaceskip=0pt\relax}
\providecommand{\BIBentryALTinterwordstretchfactor}{4}
\providecommand{\BIBentryALTinterwordspacing}{\spaceskip=\fontdimen2\font plus
\BIBentryALTinterwordstretchfactor\fontdimen3\font minus
  \fontdimen4\font\relax}
\providecommand{\BIBforeignlanguage}[2]{{%
\expandafter\ifx\csname l@#1\endcsname\relax
\typeout{** WARNING: IEEEtran.bst: No hyphenation pattern has been}%
\typeout{** loaded for the language `#1'. Using the pattern for}%
\typeout{** the default language instead.}%
\else
\language=\csname l@#1\endcsname
\fi
#2}}
\providecommand{\BIBdecl}{\relax}
\BIBdecl

\bibitem{andriulli2008multiplicative}
F.~P. Andriulli, K.~Cools, H.~Bagci, F.~Olyslager, A.~Buffa, S.~Christiansen,
  and E.~Michielssen, ``A multiplicative calderon preconditioner for the
  electric field integral equation,'' \emph{IEEE Transactions on Antennas and
  Propagation}, vol.~56, no.~8, pp. 2398--2412, 2008.

\bibitem{andriulli2010solving}
F.~P. Andriulli, A.~Tabacco, and G.~Vecchi, ``Solving the efie at low
  frequencies with a conditioning that grows only logarithmically with the
  number of unknowns,'' \emph{IEEE Transactions on Antennas and Propagation},
  vol.~58, no.~5, pp. 1614--1624, 2010.

\bibitem{andriulli2009hybrid}
F.~P. Andriulli, H.~Bagci, K.~Cools, E.~Michielssen, F.~Olyslager, and
  G.~Vecchi, ``An hybrid calder{\'o}n-hierarchical preconditioner for the efie
  analysis of radiation and scattering from pec bodies.'' in \emph{2009 IEEE
  Antennas and Propagation Society International Symposium}.\hskip 1em plus
  0.5em minus 0.4em\relax IEEE, 2009, pp. 1--4.

\bibitem{lee2003incomplete}
J.~Lee, J.~Zhang, and C.-C. Lu, ``Incomplete lu preconditioning for large scale
  dense complex linear systems from electromagnetic wave scattering problems,''
  \emph{Journal of Computational Physics}, vol. 185, no.~1, pp. 158--175, 2003.

\bibitem{benzi1996sparse}
M.~Benzi, C.~D. Meyer, and M.~Tuma, ``A sparse approximate inverse
  preconditioner for the conjugate gradient method,'' \emph{SIAM Journal on
  Scientific Computing}, vol.~17, no.~5, pp. 1135--1149, 1996.

\bibitem{adrian2019refinement}
S.~B. Adrian, F.~P. Andriulli, and T.~F. Eibert, ``On a refinement-free
  calder{\'o}n multiplicative preconditioner for the electric field integral
  equation,'' \emph{Journal of Computational Physics}, vol. 376, pp.
  1232--1252, 2019.

\bibitem{rao1982electromagnetic}
S.~Rao, D.~Wilton, and A.~Glisson, ``Electromagnetic scattering by surfaces of
  arbitrary shape,'' \emph{IEEE Transactions on antennas and propagation},
  vol.~30, no.~3, pp. 409--418, 1982.

\bibitem{andriulli2012loop}
F.~P. Andriulli, ``Loop-star and loop-tree decompositions: Analysis and
  efficient algorithms,'' \emph{IEEE Transactions on Antennas and Propagation},
  vol.~60, no.~5, pp. 2347--2356, 2012.

\end{thebibliography}

\end{document}